\documentclass{article} 
\usepackage{amssymb}
\usepackage{amsmath}
\begin{document} 
\input epsf.sty

\title{On the functional equation $F(A(z))=G(B(z)),$ where $A,B$ are polynomials 
and $F,G$ are continuous functions 
}
 
\author{F. Pakovich}
\date{}

\maketitle

\def\be{\begin{equation}}
\def\ee{\end{equation}}
\def\bs{$\square$ \vskip 0.2cm}
\def\d{{\rm d}} 
\def\D{{\rm D}} 
\def\I{{\rm I}} 
\def\C{{\mathbb C}} 
\def\N{{\mathbb N}} 
\def\P{{\mathbb P}}
\def\Z{{\mathbb Z}}
\def\R{{\mathbb R}} 
\def\ord{{\rm ord}}
\def\ssigma{\omega}

\def\e{\eqref}
\def\phi{{\varphi}}
\def\v{{\varepsilon}} 
\def\deg{{\rm deg\,}} 
\def\Det{{\rm Det}}
\def\dim{{\rm dim\,}} 
\def\Ker{{\rm Ker\,}} 
\def\Gal{{\rm Gal\,}}
\def\St{{\rm St\,}} 
\def\exp{{\rm exp\,}} 
\def\cos{{\rm cos\,}} 
\def\diag{{\rm diag\,}} 
\def\GCD{{\rm GCD }}
\def\LCM{{\rm LCM }}
\def\mod{{\rm mod\ }}

\def\bp{\begin{proposition}}
\def\ep{\end{proposition}}
\newtheorem{zzz}{Theorem}
\newtheorem{yyy}{Corollary}
\def\bt{\begin{theorem}}
\def\et{\end{theorem}}
\def\be{\begin{equation}}
\def\bee{\begin{equation*}}
\def\la{\label}
\def\l{\lambda}
\def\m{\mu}
\def\ee{\end{equation}}
\def\eee{\end{equation*}}
\def\bl{\begin{lemma}}
\def\el{\end{lemma}}
\def\bc{\begin{corollary}}
\def\ec{\end{corollary}}
\def\pr{\noindent{\it Proof. }}
\def\note{\noindent{\bf Note. }}
\def\bd{\begin{definition}}
\def\ed{\end{definition}}
\def\e{\eqref}

\newtheorem{theorem}{Theorem}[section]
\newtheorem{lemma}{Lemma}[section]
\newtheorem{definition}{Definition}[section]
\newtheorem{corollary}{Corollary}[section]
\newtheorem{proposition}{Proposition}[section]


\section{Introduction}
In this note we describe solutions of
the equation \be \la{00} F(A(z))=G(B(z)),\ee where $A,B$ are polynomials 
and $F,G \,:\, \C\P^1 \rightarrow \C\P^1$are non-constant continuous functions on 
the Riemann sphere. Our main result is the following theorem.  
\vskip 0.2cm
\noindent{\bf Theorem.} {\it 
Let $A,B$ be complex polynomials
and $F,G \,:\, \C\P^1 \rightarrow \C\P^1$ be non-constant continuous functions  
such that 
equality \eqref{00} holds for any $z\in \C\P^1$. Then there exist  
polynomials $C,D$ such 
\be \la{r} C(A(z))=D(B(z)). \ee
Furthermore, 
there exists a continuous function $H \,:\, \C\P^1 \rightarrow \C\P^1$ 
such that 
\be \la{co} F(z)=H(C(z)), \ \ \ \ G(z)=H(D(z)).\ee
}

Note that since all polynomial solutions of equation \eqref{r} are described by Ritt's 
theory of factorisation of polynomials (see \cite{ri}, \cite{sch}) 
the theorem above  
provides essentially complete solution of the problem. Note also that 
if the functions $F,G$ are rational then the function $H$ is also rational (see Remark below).

The idea of our approach is to use a recent result of \cite{p} which describes the
collections $A,B,K_1,K_2,$ where $A,B$ polynomials
and $K_1,K_2$
are infinite compact subsets of $\C$
such that the condition
\be \la{1} A^{-1}\{K_1\}=B^{-1}\{K_2\}\ee holds. It was shown in 
\cite{p} that \eqref{1} implies that there exist polynomials 
$C,D$ and a compact set $K\subset \C$ such that \eqref{r} holds
and 
$$ K_1=C^{-1}\{K\}, \ \ \ \ K_2=D^{-1}\{K\}.$$

The connection of \eqref{1} and \eqref{00} is clear: if equality \eqref{00} holds then 
for any set $K\subset \C\P^1$ equality \eqref{1} holds with 
\be \la{k} K_1=F^{-1}\{K\}, \ \ \ \ K_2=G^{-1}\{K\}. \ee
Therefore, if $F,G $
are any functions $\C\P^1 \rightarrow \C\P^1$ or $\C \rightarrow \C\P^1$
such that
there exists
a set $K\subset \C\P^1$
for which $F^{-1}\{K\}$ and $G^{-1}\{K\}$ are infinite compact subsets of $\C$
then the result of \cite{p} permits to conclude that
equality \eqref{00} for some polynomials $A,B$  
implies that there exist polynomials $C,D$ such that equality \eqref{r} holds.  

Note however that the condition above does not hold for all interesting classes of functions. For instance,
for any meromorphic transcendental function on $\C$
the preimage of any non-exceptional value is infinite and therefore unbounded
and equation \eqref{00} where $F,G$ are function meromorphic on $\C$
in general does not imply that \eqref{r} 
holds (see \cite{l}).

\section{Proof of the theorem}
First of all observe that, since $F,G$ are continuous and $\C\P^1$ is a connected compact set, 
the set $R=F(\C\P^1)=G(\C\P^1)$ is a
connected compact set. Let now $t$ be any point of $R$ distinct from $s=F(\infty)=G(\infty)$
and $C$ be a disk with center at $t$ which  
does not contain $s.$ Set $K=R\cap C.$

Since $R$ is connected and contains more than one point the set $K$ is infinite. Besides, in view of compactness of $R$ the set $K$ is closed. Finally, any of sets 
$K_1=F^{-1}\{K\},$ $K_2=G^{-1}\{K\}$ is bounded. Indeed, if say a sequence 
$x_n\in K_1$ converges to the infinity then, since $K$ is closed, the continuity of $F$ implies that $F(\infty)\in K$ in contradiction with the initial assumption. 

It follows that $K_1, K_2$ 
are infinite compact subsets of $\C$ 
for which equality \eqref{1} holds. Set $a=\deg A(z), b=\deg B(z)$ and suppose 
without loss of generality that $a\leq b.$
By Theorem 1 of \cite{p} equality \eqref{1} implies that if $a$ divides $b$
then there exists a polynomial $C(z)$ such that $B(z)=C(A(z))$ while if $a$ does not divide $b$  
then there exist polynomials $C,D$ such that equality \eqref{r} holds. Furthermore, in the 
last case there exist polynomial $W,$ $\deg W=w=\GCD(a,b)$ and a linear functions $\sigma$ such that 
$$ A(z)=\tilde A(W(z)), \ \ \ \ B(z)=\tilde B(W(z)),\ \ \ \   
$$ where either
$$ C(z)=z^cR^{a/w}(z) \circ \sigma^{-1}, \ \ \
\tilde A(z)=\sigma \circ z^{a/w} , $$
\be \la{us1} \ \ \  \ \ \    D(z)=z^{a/w} \circ \sigma^{-1}, \ \ \ \ \ \ \ \ \ \ \tilde B(z)=\sigma \circ
z^cR(z^{a/w})\ee 
for some polynomial $R(z)$ and $c\geq 0,$ 
or
$$ C(z)=T_{b/w}(z) \circ \sigma^{-1}, \ \ \
\tilde A(z)=\sigma \circ T_{a/w}(z), $$ \be \la{us2}
D(z)=T_{a/w}(z) \circ \sigma^{-1}, \ \ \ \tilde B(z)=\sigma\circ
T_{b/w}(z)\ee
for the Chebyshev polynomials $T_{a/w}(z), T_{b/w}(z).$ 

If $a$ divides $b$ then we  
have: 
$$F\circ A=G \circ B=G\circ C \circ A.$$ Therefore,
$F(z)=G(C(z))$ and equalities \eqref{r}, \eqref{co} hold with $D(z)=z.$ So, in the following we will 
assume that $a$ does not divide $b.$

Set $U=F \circ \sigma, V=G \circ \sigma.$ Then either 
equality \be \la{l1} U \circ z^{d_1}= V \circ z^{e}R(z^{d_1})\ee or 
equality \be \la{l2} U \circ T_{d_1}(z)=V \circ T_{d_2}(z)\ee
holds with $d_1=a/w,$ $d_2=\deg\, z^{e}R(z^{d_1})=b/w.$

Since $d_1,d_2$ are coprime the theorem follows now from the following lemma.
\vskip 0.2cm
\noindent{\bf Lemma} {\it   
Let 
$U,V\,:\, \C\P^1 \rightarrow \C\P^1$ be functions 
such that equality \eqref{l1} (resp. \eqref{l2}) holds 
with $d_1$ and  
$d_2$ coprime.
Then there exists a function $H\,:\, \C\P^1 \rightarrow \C\P^1$ such that the equalities
\be \la{l3}  U(z)=H\circ z^{e}R^{d_1}(z), \ \ \ \ V(z)=H\circ z^{d_1}  \ee
(resp.  
\be \la{l4} U(z)=H\circ T_{d_2}(z), \ \ \ 
V(z)=H\circ T_{d_1}(z)\ \ {\it )} \ee
hold. Furthermore, if the functions $U,V$ are continuous
then the function $H$ is also continuous. 
}
\vskip 0.2cm
\noindent{\it Proof of the lemma.} We use the following observation (cf. \cite{e2}). Let $X$ be an arbitrary set and $f,g\ :\ X\rightarrow X$ be two functions. Then $f=h(g)$ for some function $h\ :\ X\rightarrow X$ if and only if for any two points $x,y\in X$ such that $g(x)=g(y)$ the equality $f(x)=f(y)$ holds. Indeed,   
in this case we can define $h$ by the formula $h(z)=f(g^{-1}(z))$.
Furthermore, if $X=\C\P^1$ and $f,g$ are continuous then 
it is clear that $h$ is also continuous. Note also that if 
$f,g$ are rational functions on $\C\P^1$ then the function 
$H(z)$ is also rational.

Consider first the case when equality \eqref{l1} holds. Suppose that for some $z_1,z_2\in \C$  
we have: $$ z_1^{d_1}=z_2^{d_1} $$ 
and let $\theta\in \C$ be a point such that $\theta^{e}R(\theta^{d_1})=z_1$
Since $d_1$ and $d_2$ are coprime, the numbers
$e$ and $d_1$ also coprime. Therefore, there
exists a $d_1$-root of unity $\varepsilon$ such that $
(\varepsilon \theta)^{e}R((\varepsilon \theta)^{d_1})=z_2.$ 

Hence,
$$V(z_1)=V(\theta^{e}R(\theta^{d_1}))=U(\theta^{d_1})=U((\varepsilon \theta)^{d_1})
=V((\varepsilon \theta)^{e}R((\varepsilon \theta)^{d_1}))=V(z_2)$$
and therefore $V=H(z^{d_1})$ for some continuous
function $H.$ Furthermore, we have:
$$U\circ z^{d_1}=V\circ z^eR(z^{d_1})=H\circ
z^{d_1} \circ z^eR(z^{d_1})=H\circ z^eR^{d_1}(z)\circ z^{d_1}. $$
Therefore, $U=H\circ z^eR^{d_1}(z).$

Consider now the case when equality \eqref{l2} holds. Let $z_1,z_2 \in \C$ be points such that \be \la{ch} T_{d_1}(z_1)=T_{d_1}(z_2)\ee 
and let $\phi\in \C$ be a point
such that $\cos\, \phi =z_1.$ Set $t_1=\cos( \phi/d_2).$ Then, since $T_n(\cos\, z)=\cos nz$, the equality $T_{d_2}(t_1)=z_1$ holds.

It follows from \eqref{ch} that 
$z_2$ has the form $z_2=\cos(\phi+2\pi k/d_1)$ for some $k=1, ..., d_1-1.$ 
Furthermore, since $d_1$ and $d_2$ are coprime, there exists a number $l$ such that 
$d_2l\equiv k\ \mod d_1.$
Therefore, for 
$t_2=\cos(\phi/d_2+2\pi l/d_1)$ the equality 
$T_{d_2}(t_2)=z_2$ holds. Besides, clearly $T_{d_1}(t_2)= T_{d_1}(t_1).$

Now we have:
$$V(z_1)=V(T_{d_2}(t_1))=U(T_{d_1}(t_1))=U(T_{d_1}(t_2))
=V(T_{d_2}(t_2))=V(z_2)$$ and therefore
$V=H(T_{d_1})$ for some continuous 
function $H.$ Furthermore, we have:
$$U\circ T_{d_1}=V\circ T_{d_2}=H\circ
T_{d_1} \circ T_{d_2}=H\circ T_{d_2}\circ T_{d_1}. $$
Therefore, $U=H\circ T_{d_2}.$

\vskip 0.4cm

\vskip 0.2cm
\noindent {\bf Remark.} As it was remarked in the proof of the lemma if functions $F,G$ in \eqref{00} are rational then the function $H$ is also rational 
and it is clear that an appropriate modification of the theorem is hold for any functions $F,G \,:\, \C\P^1 \rightarrow \C\P^1$ or $\C \rightarrow \C\P^1$ such that there exists
a set $K\subset \C\P^1$
for which $F^{-1}\{K\}$ and $G^{-1}\{K\}$ are infinite compact subsets of $\C.$


Note however that if $F,G$ are rational then the theorem can be establish 
much easier (see also \cite{l}, Th. 2 for some more general approach to equation \eqref{00} with rational $F,G$). Indeed,
if $F=F_1/F_2$, where polynomials $F_1,F_2$ have no common roots then polynomials  $F_1(A),F_2(A)$ also have no common roots.
Similarly, if $G=G_1/G_2,$ where polynomials $G_1,G_2$ have no common roots then polynomials $G_1(B),G_2(B)$ have no common roots. Therefore, if equality \eqref{00} holds then there exists $c\in \C$ such that equality \eqref{r} holds
with $C(z)=F_1(z)$ and $D(z)=cG_1(z).$ Furthermore, it follows from the Ritt theorem
(\cite{ri}, \cite{sch}) that if $C,D$ are polynomials of minimal degrees satisfying \eqref{r}
then $\deg C=b/w,$ $\deg D=a/w.$ 
This implies that 
$$\C(C(z))\cap\C(D(z))=\C(W(z)),$$
where $$W(z)=C(A(z))=D(B(z)).$$
Now the L\"uroth theorem implies easily that equalities \eqref{co} hold.

\vskip 0.2cm
\noindent{\bf Acknowledgments.} I am grateful to
A. Eremenko for interesting discussions.

\end{document}